\documentclass[review]{elsarticle}
\usepackage{lineno,hyperref}

\journal{International Journal of Number Theory}
\usepackage{amssymb}
\usepackage{setspace}
\usepackage{fullpage}
\usepackage{amsfonts}
\usepackage{amsmath}
\usepackage{comment}
\usepackage{amsthm}
\usepackage{amssymb,eucal,mathrsfs,latexsym,xy,makeidx,verbatim}
\theoremstyle{plain}  
\newcommand{\ndiv}{\hspace{-4pt}\not|\hspace{2pt}}
\newcommand{\Mod}[1]{\ (\mathrm{mod}\ #1)}

\newtheorem{theorem}{Theorem}[section]

\newtheorem{lemma}[theorem]{Lemma}

\newtheorem{proposition}[theorem]{Proposition}

\newtheorem{rem}[theorem]{Remark}
\theoremstyle{definition}

\newtheorem{example}[theorem]{Example}

\begin{document}
\begin{frontmatter}
\title{Gauss-Dickson Codes}

\author[1]{S. A. Katre}
\address[1]{Department of Mathematics, SPPU, Pune-411001, India.}
 \ead{sakatre@gmail.com}

\author[2]{Vikas S. Jadhav \corref{cor1}}
 \ead{svikasjadhav@gmail.com}
\cortext[cor1]{Corresponding author: svikasjadhav@gmail.com}
\begin{abstract}	
Let $l$ be an odd prime. For primes, $p\equiv 1 \Mod{l}$, Gauss $(l=3)$ and Dickson $(l=5)$ 
considered the Diophantine systems in terms of which cyclotomic numbers of order 3 and 5 were obtained. The aim of this paper is to show how to obtain 1-error detecting $[2, 1, 2]$ code and 1-error correcting $[4, 2, 3]$ code in terms of the solutions of these diophantine systems in the set up of finite fields of $q=p^{\alpha}$ elements,  $p\equiv 1 \Mod{l}$, $l=3,5$.
\end{abstract}


\end{frontmatter}


\section{Introduction}
Gauss used his cyclotomic periods to obtain cyclotomic numbers of order 3 and 4, however Dickson used Jacobi sums to obtain cyclotomic numbers of order 3, 4, 5 and more  (See \cite{bruce} for details).
\vskip2mm \noindent
 \textbf{Diophantine systems of Gauss:} In $1801$, Gauss published his famous book Disquisitiones Arithmeticae in which he 
introduced cyclotomic numbers of order $3$ and $4$, as an incidental application of his theory of cyclotomy, 
which he developed to solve the longstanding problem of constructibility of regular $n$-gons. 
Gauss used Gaussian periods or Gauss sums for this work and he obtained cyclotomic numbers of 
order 3 and 4 in terms of solutions of  the diophantine systems for a prime $p$:
\begin{align*}
	p \equiv 1\Mod{3} ;\,\, 4p=L^2+27M^2,\,\, L \equiv 1\Mod{3}.\\
	p \equiv 1\Mod{4} ;\,\, p= a^2 + b^2,\,\, a \equiv 1\Mod{4}.
\end{align*}
\noindent In the set up $q=p^{\alpha}$,  $p$ prime $\equiv 1 \Mod{3}$, the system becomes $4q=L^2+27M^2,\,\, L \equiv 1\Mod{3}, p\nmid L$ (See \cite{kr}). 
For $q=p^{\alpha}$,  $p$ prime $\equiv 1 \Mod{4}$, we have $q=a^2+b^2$, $a\equiv 1 \Mod{4}$, $ p\nmid a$.
\vskip2mm \noindent
\textbf{Dickson's Work:} Around $1935$, L. E. Dickson extended the results of Gauss on cyclotomic numbers and  used them for his work on Waring's problem.
\vskip2mm \noindent
Dickson used Jacobi sums to study cyclotomic numbers. Using properties of Jacobi sums of order $5$ he obtained a diophantine system
for primes $p \equiv 1 \Mod 5$:
\begin{align*}
	16p= X^2 +50U^2+50V^2+125W^2,\hspace{4cm} (2) \\
	XW=V^2-4UV-U^2, X\equiv 1 \Mod 5. \hspace{4cm}
\end{align*}
This system has 4 solutions $ (X,U,V,W)$.
\vskip2mm \noindent
\textbf{Parnami, Agrawal, Rajwade \cite{par}:} These authors generalised Dickson's diophantine system for $q=p^\alpha$, $p\equiv 1 \Mod 5$.\\
\textbf{Katre, Rajwade \cite{kr1}, \cite{kr2}:} These authors simplified the conditions of Parnami, Agrawal and Rajwade and reduced the 
$(\alpha +1)^2$ solutions of the Dickson's system for  $q=p^\alpha$, $p\equiv 1 \Mod 5$ to $4$ solutions by the additional condition
$p \not | (X^2-125W^2)$. Thus we get the diophantine system 
\begin{align*}
16q= X^2 +50U^2+50V^2+125W^2, \hspace{4cm} (3)\\
XW=V^2-4UV-U^2, X\equiv 1 \Mod 5,\,\, p \not | (X^2-125W^2).
\end{align*}
They also showed that Dickson's formulae for cyclotomic numbers of order $5$ related to a 
given generator $\gamma$ work by choosing a unique solution of the system $(2)$ by a new condition: 
$\gamma^{(q-1)/5}\equiv \frac{A-10B}{A+10B} \Mod p$. 
This resolved the ambiguity in the determination of cyclotomic numbers of order $5$. \\
We shall now see how to get MDS codes from the generalised Gauss and Dickson systems over $\mathbb F_q$.
\section{Basic results} 

\subsection{Arithmetic characterization of Jacobi sums} 
Our development of MDS codes is based on an arithmetic characterisation of Jacobi sums of order $l$ given by Katre-Rajwade \cite{kr}.
We state important properties of Jacobi sums of order $l$ and their arithmetic characterisation.\\
Let $l$ be an odd  prime. Let $\mathbb Z$ be the ring of rational integers,  
$\mathbb Q$ the field of rational numbers, $\zeta_l =$ exp$(2\pi i/l)$. The ring of algebraic integers 
in the cyclotomic field $\mathbb Q(\zeta_l)$ is $\mathbb Z[\zeta_l]$ and it is a Dedekind domain.
 The only roots of unity in $\mathbb Z[\zeta_l]$ are $\pm \zeta_l^i, ~0 \leq i \leq l-1$. The units of $\mathbb Z[\zeta_l]$ are 
$$ \pm \zeta_l^i \displaystyle \prod_{a}\left(\zeta_l^{\frac{1-a}{2}}\frac{1-\zeta_l^a}{1-\zeta_l}\right)^{j_a}, 
~1 < a \leq \frac{l-1}{2}, ~(a,p)=1,  ~i,j_a \in \mathbb Z, ~0 \leq i \leq l-1 \,\,(\mbox{see} \,\,\cite{lw}).$$ 
$(1- \zeta_l)$ is a prime ideal in $\mathbb Z[\zeta_l]$ and $(l) = (1- \zeta_l)^{l-1}$ as ideals. 
\vskip2mm \noindent
Let $p$ be a prime $\equiv 1 \Mod{l}$. Then $p$ splits completely in $\mathbb Q[\zeta_l]$. 
The Galois group Gal$(\mathbb Q(\zeta_l)/ \mathbb Q)$ is cyclic and it consists of the automorphisms 
$\sigma_i ~(1 \leq i \leq l-1)$, defined by $\sigma_i(\zeta_l) = \zeta_l^i$.
Then $(p) =\displaystyle \prod_{i=1}^{l-1} \mathcal P^{\sigma_i}$ where 
$\mathcal P$ be a prime ideal factor of $p$ in  $\mathbb Z[\zeta_l]$. 
Also, $N(\mathcal P)= p$, where $N(\mathcal P)$ denotes the norm of the ideal $\mathcal P$. 
Hence $(\mathbb Z[\zeta_l]/\mathcal P)^* \cong (\mathbb Z / p \mathbb Z)^*$ is cyclic of order $p-1$. 
Let $\gamma$ be a generator of $\mathbb F_q^*, q=p^{\alpha}, \alpha \geq 1$. 
Let $b = \gamma^{\frac{q-1}{l}}$. As $p \equiv 1 \Mod{l},\,\, b \in \mathbb F_p^*$ and 
we can take it to be an integer $\Mod p$. We have,
\begin{equation*}
 0 \equiv \gamma^{q-1} - 1 \Mod{\mathcal P} \\
\end{equation*}
\begin{equation*}
 ~~~~~~~~~\equiv \displaystyle \prod_{i=0}^{l-1}(\gamma^{\frac{q-1}{l}} - \zeta_l^i) \Mod{\mathcal P}.
\end{equation*}
Since $l^{\rm{th}}$ roots of unity are distinct $\Mod{\mathcal P}$, we get $\gamma^{\frac{q-1}{l}} \equiv \zeta_l^k \Mod{\mathcal P} $
for exactly one $k, ~1 \leq k \leq l - 1$. For the given prime factor $\mathcal P$, we can choose a primitive
root $\gamma$ so that $k = 1$. On the other hand if we first choose a generator $\gamma$ of $\mathbb F_q^*$, 
there is a unique prime divisor $\mathcal P$ of $(p)$ in $\mathbb Z[\zeta_l]$, such that
\begin{equation}\label{eqn521}
 \gamma^{\frac{q-1}{l}} \equiv \zeta ~(\mbox{mod} ~\mathcal P). 
\end{equation}
We thus assume that $\gamma$ and $\mathcal P$ are related by (\ref{eqn521}). We define the character $\chi_l$ on 
$\mathbb F_q^{*}$ by  $\chi_l(\gamma) = \zeta_l$. Define
\begin{equation*}
J(i,j)_l = \displaystyle \sum_{-1 \neq v \in \mathbb F_p^*}\chi_l^i(v)\chi_l^j(v+1).
\end{equation*}
$ J(i,j)_l$ is called a Jacobi sum of order $l$.
To know more about Jacobi sums one may refer to \cite{bruce}.
\vskip2mm \noindent
Let $\psi = J(1,1)_l $ and $\psi_i = \sigma_i(\psi), ~ 1 \leq i \leq l-1$. 
Then $\psi$ satisfies the following properties (see \cite{par}):
\begin{lemma} \label{mjlm1}
 $\psi \bar{\psi} = q$.
\end{lemma}
\begin{lemma}\label{mjlm2}
 $\psi \equiv -1 \Mod{(1- \zeta_l)^2}$.
\end{lemma}
%
%
\begin{lemma}(\cite{kjt}, \cite{par})\label{mjlm3}
 $GCD(\psi_1, \cdots, \psi_{\frac{l-1}{2}})$ is the $\alpha^{\text{th}}$ power of a prime ideal $\mathcal P$ of $\mathbb Z[\zeta_l]$. 
If $J(1, 1)_l$ is defined in terms of $\gamma$, then $\mathcal P$ coincides with the prime ideal described by (\ref{eqn521}).
Moreover 
\begin{equation*}
 (\psi)= (J(1, 1)_l)= \displaystyle \prod_{k=1}^{(l-1)/2} (\mathcal P^{\displaystyle\sigma_{k^{-1}}})^\alpha
\end{equation*}
where $k^{-1}$ is taken (mod $l$). 
\end{lemma}
\noindent Similar factorisation can also given for $(J(1, n)), \,n > 1$. Lemmas \ref{mjlm1}, \ref{mjlm2} and \ref{mjlm3}
together give an algebraic characterisation of the Jacobi sum $J(1, 1)$ as an element of $\mathbb Z[\zeta_l]$.
\vskip1mm 
Let $H = \displaystyle \sum_{i=0}^{l-1}a_i(n) \zeta_l^i$ with $a_0(n) = 0, ~1 \leq n \leq l-2$.
The suffixes in $a_i(n)$ are to be considered $\Mod l$.
Parnami, Agrawal and Rajwade \cite{par} showed that for each
$n, ~1 \leq n \leq l - 2$, the Diophantine system (arithmetical conditions):
\vskip1mm
(i) $~~~~q = \displaystyle \sum_{i=1}^{l-1}a_i^2(n) -  \sum_{i=1}^{l-1}a_i(n)a_{i+1}(n)$,
\vskip1mm
(ii) $~~\displaystyle \sum_{i=1}^{l-1}a_i(n)a_{i+1}(n) = \sum_{i=1}^{l-1}a_i(n)a_{i+2}(n) = \cdots = \sum_{i=1}^{l-1}a_i(n)a_{i+l-1}(n)$,
\vskip1mm
(iii) $~~1 + \displaystyle \sum_{i=1}^{l-1}a_i(n) \equiv 0\Mod{l}$, 
\vskip1mm
(iv) $~~\displaystyle \sum_{i=1}^{l-1}i a_i(n) \equiv 0 \Mod{l}$,
\vskip1mm
(v) $~~~ p  \ndiv \prod_{\lambda((n+1)k) > k} H^{\sigma_k}$, $(\lambda(r)$ being the least non-negative remainder of  $ r \Mod l )$
\vskip2mm \noindent
has  $l - 1$  solutions, so that  $H = \displaystyle \sum_{i=0}^{l-1}a_i(n) \zeta_l^i$ 
is one of the $l-1$ field conjugates of $J(1, n)_l$ and conversely. 
 Thus if $(a_1(n), a_2(n),\cdots , a_{l-1}(n))$ is a solution of the above system then its other
solutions are $(a_{i\cdot 1}(n), a_{i\cdot 2}(n), \cdots , a_{i\cdot (l-1)}(n))$, for $2 \leq i \leq l - 1$.
Here the number of distinct solutions of $(i)-(v)$ is equal to the number of distinct conjugates of $J(1, n)_l$.
For $n=1$, all the $l-1$ conjugates of $J(1, 1)_l$ are distinct and so we get $l-1$ distinct solutions of $(i) - (v)$.
\vskip2mm \noindent
Katre and Rajwade \cite{kr} showed that for an integer $b \equiv \gamma^{\frac{q-1}{l}}( \mbox{mod} ~p)$, 
 the additional condition
\vskip1mm
(vi) $~~ p \mid \overline{H}\prod_{\lambda((n+1)k) > k}(b - \zeta_l^{\sigma_{k^{-1}}}) $, 
\vskip1mm
where $k^{-1}$ is taken $(\mbox{mod} ~l)$, determines the unique solution $H = J(1, n)_l$.\\
This resolved a longstanding ambiguity in cyclotomy for determination of Jacobi sums and cyclotomic numbers of order $l$.
\vskip2mm \noindent
If $n = 1$, (v) becomes $~~~ p \ndiv \displaystyle \prod_{k=1}^{(l-1)/2} H^{\sigma_k}$ and $(vi)$ becomes
$~~ p \mid \displaystyle \overline{H}\prod_{k=1}^{(l-1)/2}(b - \zeta_l^{\sigma_{k^{-1}}}) $.
\vskip2mm 
\section{A conjecture related to Jacobi sums and the construction of MDS codes of type
$[l-1, (l-1)/2, (l+1)/2]$}

For $q=p^{\alpha}$, $p \equiv 1 \Mod l$, $J(1, 1)$ has $l-1$ distinct conjugates in $\mathbb Z[\zeta_l]$. We fix $n=1$ for the 
discussion in \S 3. For any solution $(a_1, a_2, \cdots, a_n)$ of $(i)-(v)$, or equivalently for any $a_1, a_2, \cdots, a_{l-1}$ satisfying
$J(1, 1)= \displaystyle \sum_{i=1}^{l-1}a_i\zeta^i$ for some generator $\gamma$ of $\mathbb F_q^*$, 
the $(vi)^{\mbox{th}}$ condition corresponding to the Diophantine system $(i)-(v)$ is a system of polynomial congruences having $b$ as a solution.
This system can be considered as  a system of $l-1$ linear congruence equations in $\frac{l-1}{2}$ variables having 
$(b, b^2, \cdots, b^{\frac{l-1}{2}})$ as a solution. 
Thus the system (considered as a system of linear equations) is consistent with rank at most $\frac{l-1}{2}$. 
We represent the above system of equations in the matrix form as 
$$ DX = Y $$
where $D$ is the  $(l-1) \times \frac{(l-1)}{2}$ matrix coming from the coefficients of $b, b^2, \cdots, \displaystyle b^{\frac{l-1}{2}}$.
Let $D^t$ denote the transpose of $D$.
\vskip2mm \noindent
{\large \bf Conjecture (S. A. Katre):} 
\vskip1mm \noindent
 Any $\displaystyle \frac{l-1}{2}$ rows of $D$ are linearly independent, except possibly for finitely many primes $p$, 
and $D^t$ is a generator matrix of  an MDS code over $F_q$.
\vskip2mm \noindent
The result has been verified for $l = 3, 5$ by S. A. Katre earlier (with no exception) \cite{sak}. 
In the next sections we demonstrate this using the systems of Gauss and Dickson in the set up of $\mathbb F_q$.

\section{MDS codes of type $[2, 1, 2]$  obtained from Gauss System}
Let $q=p^{\alpha}$,   $p$ prime $\equiv 1 \Mod{l}$, $\gamma$ be a generator of $\mathbb F_q^*$ and $\zeta_l$ be a primitive $l$-th root of unity in $\mathbb C$, $l$ odd prime. We note here that for $l=3,5$, these diophantine systems of Gauss and Dickson were obtained in the set up of the finite  fields $\mathbb F_q$, by Parnami-Agarwal-Rajwade-Katre using the properties of the Jacobi sums
\begin{equation*}
	J(1,1)_l = \displaystyle \sum_{-1,0 \neq v \in \mathbb F_q^*}\chi_l(v)\chi_l(v+1).
\end{equation*}
where $\chi$ is a character on $\mathbb F_q^*$ satisfying $\chi(\gamma)= \zeta_l$. We recall that an $[n, k, d]$-code over $\mathbb F_q$ is a subspace of $\mathbb F_q^n$ of dimension $k$ and distance $d$. Such a code detects $d-1$ errors and corrects $[\frac{d-1}{2}]$. This code is called an MDS (Maximum distance separable) code if $n+1=k+d$. (See \cite{vp}, \cite{ling})
\begin{proposition}\label{mjprop2}
	A $q$-ary linear $[n,k]$-code $C$ is an MDS code if and
	only if every set of $n-k$ columns of a parity check matrix of $C$ is linearly
	independent.
\end{proposition}

\begin{proposition}\label{mjprop3}
	A $q$-ary linear $[n,k]$-code $C$ is an MDS code if and only if every set of $k$ columns of a generator matrix of $C$ is linearly independent.
\end{proposition}
In the Gauss case we get a 1-error detecting MDS code of the type $[2, 2, 1]$  and in the Dickson case we get a  1-error correcting MDS code of the type $[4, 2, 3]$  for all  $q=p^{\alpha}$,  $p\equiv 1 \Mod{l}$, $l=3,5$.
\vskip2mm \noindent
\underline{\bf The case $l =3$ } (cf. \cite{kr}, \cite{par}): 
\vskip2mm \noindent
Let $\zeta = e^{ \frac{2\pi i}{3}},$ $q=p^{\alpha}, p \equiv 1 \Mod 3$, $p$ a prime. Let for a generator 
$\gamma$ of $\mathbb F_q^*, b=\gamma^{ \frac{q-1}{3}}$. Then $b \in \mathbb F_p$,
which we take as an integer $\Mod p$. $b$ is a cube root of unity $\Mod p$. The diophantine system  $(i)-(vi)$ considered by Katre-Rajwade, (see \cite{kr}, \cite{par}), for the Jacobi sum 
$J(1, 1) = a_1 \zeta + a_2 \zeta^2 $ of order 3  takes the form\\
$i)$ \,\,$\,\,q=a_1^2+a_2^2 - a_1a_2$\\
$ii)$ \,\,no condition in this case \\
$iii)$ \,\,$1+a_1+a_2\equiv 0 \Mod 3$\\
$iv)$ \,\,$a_1+2a_2\equiv 0 \Mod 3 $ i.e. $a_1 -a_2 \equiv 0 \Mod 3$\\
$v)$\,\,$ \,\, p \not | (a_1 \zeta + a_2 \zeta^2) $\\
together with a $(vi)^{\mbox{th}}$ condition:\\
$vi)$ $a_2 b + a_1\equiv 0 \Mod p,$\\
\hspace*{5mm} $ a_1 b + a_1 - a_2 \equiv 0 \Mod{p}$. 
\vskip1mm \noindent
We use the transformations (between the solutions of diophantine systems) 
$$a_1 = \frac{-L+3M}{2}, \,\,a_2 = \frac{-L-3M}{2} $$
 with the inverse transformations 
$$L = -(a_1 + a_2),\,\, M = \frac{a_1 -a_2}{3}.$$
 Then the system $(i)-(v)$ takes the form  of the (generalised) Gauss system:
 $$4q = L^2 + 27M^2, \,\, p \not| \,\,L,\,\, L \equiv 1 \Mod p.$$
This determines $L$ uniquely and $M$ up to sign. For $q=p$, the condition $p\nmid L$ is automatically satisfied. Then 
$$J(1, 1)= \frac{L+3M}{2}+3M\zeta = \frac{-L+3M}{2}\zeta + \frac{-L-3M}{2}\zeta^2,$$
where the sign of $M$ is determined by condition $(vi)$, which takes the form 
\begin{eqnarray}\label{eq2}
b \left(\frac{L+3M}{2}\right)-\frac{L-3M}{2} \equiv 0 \Mod p,\\
b \left(\frac{L-3M}{2}\right)-3 M\equiv 0 \Mod p.\nonumber
\end{eqnarray}
Suppose $\frac{L-3M}{2} \equiv 0 \Mod p$. 
Then by the first congruence equation we get $\frac{L+3M}{2} \equiv 0 \Mod p$. By adding we get $L\equiv{0 \Mod p}$, a contradiction as $p\not| L$. Hence $\frac{L-3M}{2}$ and so  $\frac{L+3M}{2}$ are nonzero $\Mod p$. 
From the second equation, $M \not \equiv 0 \Mod p$.
Thus each of these congruence equations has nonzero coefficients $\Mod p$  and  is linearly independent $\Mod p$.
However the two equations are linearly dependent $\Mod p$ as the system is consistent, which can also be seen using $4q = L^2 + 27M^2$.
Moreover  each of these equations and the condition $(vi)$ is equivalent to $b \equiv \frac{L-3M}{L+3M} \Mod p$, 
giving that $ \frac{L-3M}{L+3M}$ is a cube root of unity $\Mod p$. 
Let $D$ be a $2 \times 1$ column matrix with entries as coefficients of $b$ in the two congruences (\ref{eq2}), so that  $D^t = [\frac{L+3M}{2}, \frac{L-3M}{2}]$. Thus the entries in the  columns of $D^t$ are nonzero $\Mod p$.
So each column of $D^t$ is linearly independent $\Mod p$.   
Hence by proposition \ref{mjprop3}, $D^t = [\frac{L+3M}{2}, \frac{L-3M}{2}]$ is  a generator matrix of a $[2, 1, 2]$-MDS code over $\mathbb F_q$.
We call this code as a Gauss-Code. It is a $1$-error detecting MDS code. Although any row of length 2 with nonzero entries $\Mod p$ determines such a code, considering the historical importance of the Gauss system, we have illustrated how the method works beginning with $l = 3$.
\section{MDS codes of type  $[4, 2, 3]$ obtained from Dickson System}
\vskip2mm \noindent
\underline{\bf The case $l =5$ } (cf. \cite{kr}, \cite{kr1},  \cite{par}): 
\vskip2mm \noindent
Let $\zeta = e^{ \frac{2\pi i}{5}}$, $q=p^{\alpha},\,\, p \equiv 1 \Mod 5$, $p$ a prime. Let for a generator 
$\gamma$ of $\mathbb F_q^*,\,\, b=\gamma^{ \frac{q-1}{5}}$. Then $b \in \mathbb F_p$,
which we take as an integer $\Mod p$. The diophantine system $(i)-(vi)$ in \cite{kr} (See also \cite{kr1}, \cite{par}) for the Jacobi sum 
$J(1, 1) = a_1 \zeta + a_2 \zeta^2 + a_3 \zeta^3 + a_4 \zeta^4$ of order 5  takes the form:\\
$i)$ $q=a_1^2+a_2^2+a_3^2+a_4^2-\frac{1}{2} (a_1a_2+a_2a_3+a_3a_4+a_1a_3+a_2a_4+a_1a_4)$\\
$ii)$ $a_1a_2+a_2a_3+a_3a_4=a_1a_3+a_2a_4+a_1a_4$\\
$iii)$ $1+a_1+a_2+a_3+a_4\equiv 0 \Mod 5$\\
$iv)$ $a_1+2a_2+3a_3+4a_4\equiv 0 \Mod 5$\\
$v)$ $p\ndiv ~~{\rm gcd} (a_2^2+a_1a_4-a_1a_2-a_2a_4-a_1a_3, \cdots, \cdots, \cdots)$\\
together with a $(vi)^{\mbox{th}}$ condition:\\
$vi)$ $b^2a_4+b(a_1-a_2+a_3)+(a_3-a_4)\equiv 0 \Mod p,$ \hskip2.4cm \,\,$(I)$\\
\hspace*{5mm} $b^2a_3+b(a_3-a_4)+(a_2-a_4)\equiv 0 \Mod p$, \hskip3cm \,\,$(II)$\\
\hspace*{5mm} $b^2a_2+ba_1+(a_1-a_4)\equiv 0 \Mod p,$ \hskip4cm \,\,$(III)$\\
\hspace*{5mm} $b^2a_1+b(a_1-a_2+a_3-a_4)-a_4\equiv 0 \Mod p$. \hskip2.5cm \,\,$(IV)$
\vskip1mm \noindent 
We use the bijective transformations from solutions of $(i)-(v)$ to the solutions of $(2)$:
\begin{align*}
a_1 = \frac{1}{4}(-X+2U+4V+5W),\,\,\,\, & \,\,\,\,a_2= \frac{1}{4}(-X+4U-2V-5W),\hspace{3cm} (*)\\
a_3 = \frac{1}{4}(-X-4U+2V-5W), \,\,\,\, & \,\,\,\, a_4= \frac{1}{4}(-X-2U-4V+5W) \notag
\end{align*}
with the inverse transformations
\begin{align*}
X = -(a_1+a_2+a_3+a_4),\,\,\,\,\,\, & \,\,\,\,U=\frac{1}{5}(a_1+2a_2-2a_3-a_4), \hspace{4cm} \\
V = \frac{1}{5}(2a_1-a_2+a_3-2a_4),\,\,\,\, & \,\,\,\,W=\frac{1}{5}(a_1 - a_2-a_3+a_4). \notag
\end{align*} 
Then the system $(i)-(v)$ takes the form  of the Dickson-Katre-Rajwade system (see \cite{kr1})\\
$(1)$ $16q = X^2 + 50U^2 + 50V^2 + 125W^2$\\
$(2)$ $XW = V^2 -4UV -U^2$\\
$(3)$ $X \equiv 1 \Mod 5$\\
$(4)$ $p \ndiv X^2 - 125W^2$ (the rejection condition). \\
Here $(1)-(3)$ are given by Dickson for $q=p$ case and $(4)$ is added by Katre-Rajwade for a general $q=p^{\alpha}$. 
$(4)$ is automatically satisfied when $q=p$.\\
It has been further shown by Katre-Rajwade that with the next condition $(5)$, the system $(i)-(vi)$ is equivalent to the system $(1)-(5)$.\\
(5) $\gamma^{\frac{q-1}{5}} \equiv \displaystyle \frac{A-10B}{A+10B} \Mod p$, where 
$A= X^2 - 125W^2, B= 2XU-XV-25VW$.\\
Note that here for a solution of $(1)-(4)$, both $A-10B$ and $A+10B$ are nonzero $\Mod p$.
For convenience, we shall use $(vi)^{th}$ condition in terms of the $a_i's$. 
This is amenable to generalisation later for finding MDS codes of higher order $l$ of length $l-1$.
\begin{rem} \label{r1}
Let $(a_1, a_2, a_3, a_4)$ be any solution of the 
 system $(i)-(v)$. In other words $\displaystyle \sum_{i=1}^4 a_i \zeta^i$ is a conjugate of $J(1, 1)$ for a 
generator $\gamma$ of $\mathbb F_q^*$. Hence for a solution
$(a_1, a_2, a_3, a_4)$ of $(i)-(v)$, there are 3 more solutions of $(i)-(v)$ and they are 
$(a_{1\cdot i}, a_{2\cdot i}, a_{3\cdot i}, a_{4 \cdot i}), \,\,2 \leq i \leq 4$, where the suffixes are modulo 5.
Thus $(a_2, a_4, a_1, a_3)$, $(a_3, a_1, a_4, a_2)$ and $(a_4, a_3, a_2, a_1)$ are also solutions of $(i)-(v)$.
Hence properties obtained for $(a_1, a_2, a_3, a_4)$ hold for these solutions too.
\end{rem}
\noindent The $4 \times 2$ matrix $D$ in $\S 3$ coming from the coefficients of $b, b^2$ in the 4 equations $(I)-(IV)$ is
\begin{align*}\nonumber
  D &=\left(
\begin{array}{cc}
 a_1 -a_2 + a_3  & a_4 \\
 a_3 - a_4 & a_3 \\
a_1 & a_2 \\
  a_1 - a_2 + a_3 - a_4 & a_1  \\
\end{array}
\right).
 \end{align*}
Note that for any solution  $(a_1, a_2, a_3, a_4)$ of $(i)-(v)$, $b$ is a unique common solution of $(vi)$ given by 
$b \equiv \gamma^{\frac{q-1}{5}} \Mod p$, for a suitable generator $\gamma$ of $\mathbb F_q^*$ and thus $b$ is a solution of any two equations.\\
It is not straight forward to check that any 2 rows of $D$ are linearly independent $\Mod p$. For this we use the historical 
system of Dickson in the generalised form for the calculations and  proceed as follows:\\
Denote the rows of $D$ by $R_1, R_2, R_3, R_4$.
\vskip2mm \noindent $(a)$  We first show that $R_1$ and $R_2$ are linearly independent $\pmod{p}$.\\
Write the  equations $(I)$ and $(II)$ in the form.
$$b(a_1-a_2+a_3)+b^2a_4+(a_3-a_4)\equiv 0 \Mod p,\hskip2cm \,\,(I)$$
$$b(a_3-a_4)+b^2a_3+(a_2-a_4)\equiv 0 \Mod p. \hskip2.5cm \,\,(II)$$
By, Cramer's rule, we get $D_1 b = N_1$, where 
$$D_1 =\begin{vmatrix}
a_1-a_2+a_3 & a_4\\
a_3-a_4 & a_3
\end{vmatrix} \,\,\text{and}\,\,N_1=\begin{vmatrix}
-(a_3-a_4) & a_4\\
-(a_2-a_4) & a_3
\end{vmatrix}. $$
Using transformations $(*)$, we get
\begin{align*}
16N_1 &=16[-a_3^2+a_4(a_2+a_3-a_4)]\\
&=-(-X-4U+2V-5W)^2+(-X-2U-4V+5W)(-X+2U+4V-15W)\Mod p\\
&=-20U^2-20V^2-100W^2+100VW-8XU+4XV \Mod p\\
&= \frac{2}{5}(-50U^2-50V^2-250W^2+250VW-20XU+10XV) \Mod p\\
&= \frac{2}{5}(X^2-125W^2+250VW-20XU+10XV)\Mod p\\
&= \frac{2}{5}(A-10B)\Mod p
\end{align*}
As $A-10B \not \equiv 0 \Mod p$, we get $N_1 \not \equiv 0 \Mod p$. Also $D_1b = N_1$, so $D_1 \not \equiv 0 \Mod p$.
Hence the first two rows of the matrix $D$ are linearly independent $\Mod p$:
\vskip2mm \noindent $(b)$ As before, form $D_2$ as the determinant of the $2\times 2$ matrix of coefficients of $b$ and 
$b^2$ in equations $(II)$ and $(III)$, and form $N_2$ as the determinant of the $2\times 2$ matrix obtained from  constants 
and coefficients of $b^2$ in equations $(II)$ and $(III)$.
Thus we get by Cramer's rule, $D_2 b = N_2$, where 
$$D_2 =\begin{vmatrix}
(a_3-a_4) & a_3\\
a_1 & a_2
\end{vmatrix} \,\,\text{and}\,\,N_2=\begin{vmatrix}
-(a_2-a_4) & a_3\\
-(a_1-a_4) & a_2
\end{vmatrix}. $$
By transformations $(*)$, we get
\begin{align*}
16D_2 &=16[a_2(a_3-a_4)-a_1a_3]\\
&=(-X+4U-2V-5W)(-2U+6V-10W)\\
& +(-X+2U+4V+5W)(X+4U-2V+5W)\Mod p\\
&=40U^2+40V^2+200W^2 \Mod p, \,\,\text{using}\,\, (1)\,\,\text{and}\,\,(2)\,\, \text{of Dickson's system}.
\end{align*}
We have $16q= X^2 +50U^2+50V^2+125W^2 = X^2 -125W^2+50(U^2+V^2+5W^2)$. Since $p \not| X^2 -125W^2 $, 
it follows that $U^2+V^2+5W^2 \not \equiv 0\Mod p$. Hence $D_2 \not \equiv 0\Mod p$ and so $N_2 \not \equiv 0\Mod p$.
 Thus $R_2$ and $R_3$ are linearly independent $\Mod p$.  
\vskip2mm \noindent $(c)$ Proceeding as in $(a)$ and $(b)$, from equations $(I)$ and $(III)$, we get $D_3 b = N_3$, where 
$$D_3 =\begin{vmatrix}
(a_1-a_2+a_3) & a_4\\
a_1 & a_2
\end{vmatrix} \,\,\text{and}\,\,N_3=\begin{vmatrix}
-(a_3-a_4) & a_4\\
-(a_1-a_4) & a_2
\end{vmatrix}. $$
Using $a_1a_2+a_2a_3+a_3a_4=a_1a_3+a_2a_4+a_1a_4$ (see $(ii)$), we see that 
$$D_3 = a_1a_2-a_2^2+a_2a_3-a_1a_4=a_1a_3+a_2a_4-a_3a_4-a_2^2= N_2  \not \equiv 0\Mod p.$$
Thus $R_1$ and $R_3$ are linearly independent $\Mod p$. 
\vskip2mm \noindent $(d)$ Using the  equations $(I)$ and $(IV)$, we get $D_4 b = N_4$, where 
$$D_4 =\begin{vmatrix}
(a_1-a_2+a_3) & a_4\\
(a_1-a_2+a_3-a_4) & a_1
\end{vmatrix} \,\,\text{and}\,\,N_4=\begin{vmatrix}
-(a_3-a_4) & a_4\\
a_4 & a_1
\end{vmatrix}. $$
Using $(ii)$ we get
$$D_4=(a_1-a_4)^2+a_2a_3, \,\,N_4= -a_4^2-a_1a_3+a_1a_4.$$
Here if $D_4\equiv 0 \Mod p$, then  $N_4\equiv 0 \Mod p$. Hence 
$$D_4+N_4=a_1^2+a_2a_3-a_1a_3-a_1a_4 \equiv 0 \Mod p.$$
But $-N_2=a_2^2+a_3a_4-a_1a_3-a_2a_4 \not \equiv 0 \Mod p$ for any solution $(a_1, a_2, a_3, a_4)$ of $(i)-(v)$.
Using Remark \ref{r1} and letting $a_1 \rightarrow a_3, a_2 \rightarrow a_1, a_3 \rightarrow a_4, a_4 \rightarrow a_2$, we get
$a_1^2+a_2a_4-a_3a_4-a_1a_2 \not \equiv 0 \Mod p$.\\
Using $a_1a_2+a_2a_3+a_3a_4=a_1a_3+a_2a_4+a_1a_4$, we have $a_1^2+a_2a_4-a_3a_4-a_1a_2=a_1^2+a_2a_3-a_1a_3-a_1a_4$. So
$$D_4+N_4=a_1^2+a_2a_3-a_1a_3-a_1a_4 \not \equiv 0 \Mod p.$$
Hence $D_4 \not \equiv 0 \Mod p$.
Thus $R_1$ and $R_4$ are linearly independent $\Mod p$.
\vskip2mm \noindent $(e)$ From the  equations $(II)$ and $(IV)$, we get $D_5 b = N_5$, where 
$$D_5 =\begin{vmatrix}
a_3-a_4 & a_3\\
a_1-a_2+a_3-a_4 & a_1
\end{vmatrix} \,\,\text{and}\,\,N_5=\begin{vmatrix}
-(a_2-a_4) & a_3\\
a_4 & a_1
\end{vmatrix}. $$
\begin{align*}
N_5&=a_1a_4-a_1a_2-a_3a_4\\
&=a_2a_3-a_2a_4-a_1a_3 \,\,\,\,(\text{since}\,\, a_1a_2+a_2a_3+a_3a_4=a_1a_3+a_2a_4+a_1a_4)\\
&=D_2 \not \equiv 0 \Mod p.
\end{align*}
As $D_5 b \equiv N_5 \Mod p$, we have $D_5 \not \equiv 0 \Mod p$.
Thus $R_2$ and $R_4$ are linearly independent $\Mod p$. 
\vskip2mm \noindent $(f)$ From  equations $(III)$ and $(IV)$, we get $D_6 b = N_6$, where 
$$D_6 =\begin{vmatrix}
a_1 & a_2\\
a_1-a_2+a_3-a_4 & a_1
\end{vmatrix} \,\,\text{and}\,\,N_6=\begin{vmatrix}
-(a_1-a_4) & a_3\\
a_4 & a_1
\end{vmatrix}. $$
Using $(ii)$, we have $D_6=a_1^2-a_1a_2+a_2^2-a_2a_3+a_2a_4=N_1\not \equiv 0 \Mod p. $\\
Thus $R_3$ and $R_4$ are linearly independent $\Mod p$. \\
From $(a)-(f)$, any two rows of the matrix $D$ are linearly independent. We thus get
\begin{theorem}\label{mds5}
 Let $q=p^{\alpha}$, $p \equiv 1 \Mod 5$. Any two rows of the $4 \times 2 $ matrix $D$ are linearly independent. 
The matrix $G = D^t$ is a generator matrix of an MDS code of type $[4, 2, 3]$ over $\mathbb F_q$.
\end{theorem}
\noindent We call this code as a Dickson code. It is a $1$-error correcting MDS code.
\newpage
\noindent \textbf{Decoding Using Jacobi Sums}\\
 Let $J(1, 1) = a_1 \zeta + a_2 \zeta^2 + a_3 \zeta^3 + a_4 \zeta^4$ be the Jacobi sum. The generator matrix used here is: 
\begin{align*}\nonumber
	G=D^t &=\left(
	\begin{array}{cccc}
		a_1 -a_2 + a_3  & 	a_3 - a_4 & a_1 & a_1 - a_2 + a_3 - a_4 \\
		a_4	& a_3 & a_2 &  a_1 
	\end{array}
	\right)
\end{align*}
Let $Y$ be the $2\times 2$ matrix consisting of first 2 columns of $G$, thus
\begin{align*}\nonumber
	Y &=\left(
	\begin{array}{cc}
		a_1 -a_2 + a_3  & 	a_3 - a_4  \\
		a_4	& a_3
	\end{array}
	\right).
\end{align*}
In the sequel, we shall use the determinants $D_1, D_2, D_3, D_4, D_5$ which are nonzero $\Mod p$.
Then  determinant  of $Y$ is $ a_1a_3-a_2a_3+a_3^2-a_3a_4+a_4^2 =D_1\neq 0$ and 
\begin{align*}\nonumber
	Y^{-1} &=\displaystyle\frac{1}{D_1}\left(
	\begin{array}{cc}
		a_3  & 	-a_3 +	 a_4  \\
		-a_4	& a_1 -a_2 + a_3
	\end{array}
	\right)
\end{align*} 
We get a generator matrix in the standard form:
\begin{align*}\nonumber
	G'=Y^{-1}G &= \left(\begin{array}{cc} 
		I_2 & \displaystyle\frac{1}{D_1}\left(
		\begin{array}{cc}
			a_3  & 	-a_3 +	 a_4  \\
			-a_4	& a_1 -a_2 + a_3
		\end{array}
		\right) \cdot\left(
		\begin{array}{cc}
			a_1 & a_1 - a_2 + a_3 - a_4 \\
			a_2 &  a_1 
		\end{array}
		\right)
	\end{array} \right)
\end{align*}
The parity check matrix is 
\begin{align*}\nonumber
	H &= \left(\begin{array}{cc} 
		\displaystyle\frac{1}{D_1} \cdot\left(
		\begin{array}{cc}
			a_1 & a_2 \\
			a_1 - a_2 + a_3 - a_4 &  a_1 
		\end{array}
		\right) \cdot \left(
		\begin{array}{cc}
			-a_3  & 	a_4  \\
			a_3 -a_4	& -a_1+a_2- a_3
		\end{array}
		\right) & I_2
	\end{array} \right)
\end{align*}
Syndrome of a received word $v=[a, b, c, d]$ is 
\begin{align*}\nonumber
	vH^t &= (a, b, c, d)\cdot \left(\begin{array}{c} 
		\displaystyle\frac{1}{\triangle} \cdot\left(
		\begin{array}{cc}
			-a_3  & a_3 -a_4	  \\
			a_4		& -a_1+a_2- a_3
		\end{array}
		\right) \cdot \left(
		\begin{array}{cc}
			a_1 & a_1 - a_2 + a_3 - a_4 \\
			a_2 &  a_1 
		\end{array}
		\right) \\
		I_2
	\end{array} \right)\\
    &= (a, b, c, d)\cdot \displaystyle \frac{1}{D_1}\left(
    \begin{array}{cc}
    	D_2 & D_5 \\
    	-D_3 &  -D_4 \\
    	1 &0\\
    	0& 1
    \end{array}
    \right)\\
	&=(a, b, c, d)\cdot \left(
	\begin{array}{cc}
		A_1 & A_2 \\
		A_3 &  A_4 \\
		1 &0\\
		0& 1
	\end{array}
	\right),
\end{align*}
where  $A_1 = \displaystyle \frac{D_2}{D_1}$, $A_2 = \displaystyle \frac{D_5}{D_1}$, $A_3 = \displaystyle \frac{-D_3}{D_1}$, $A_4 = \displaystyle \frac{-D_4}{D_1}$ are non-zero $\pmod{p}$, as $D_1$, $D_2$, $D_3$, $D_4$ and $D_5$ are non-zero $\pmod{p}$.
Thus $vH^t=	(aA_1 +bA_3+c \,\,\,\, aA_2 +bA_4+d)$. Since our code is a 1-error correcting code, we have the syndrome-decoding table of the form:
\begin{center}
	\begin{tabular}{c|c}
		syndrome of $v$ & error vector\\
		\hline\\
		$(e_0A_1, e_0A_2)$ & $(e_0, 0, 0, 0)$\\
		$(e_0A_3, e_0A_4)$ & $(0, e_0, 0, 0)$\\
		$(e_0, 0)$ & $(0, 0, e_0, 0)$\\
		$(0, e_0)$ & $(0, 0, 0, e_0)$
	\end{tabular}	
\end{center}
\begin{example}
Let  $ p = 61$.
\vskip1mm \noindent
We explicitly compute a generating matrix of an MDS code of type $[4, 2, 3]$ over $\mathbb F_{61}$. 
We take $\gamma = 2$ as the primitive root in $\mathbb F_{61}$, i.e. a generator of the cyclic group 
$\mathbb F_{61}^*$, and keeping the same notations as above we get
\begin{eqnarray*}
	J(1, 1)_{5} &=& a_1\zeta_{5} + a_2 \zeta_{5}^2  + a_{3} \zeta_{5}^{3}
	+ a_{4} \zeta_{5}^{4}\\
	&=& -6 \zeta_{5}^2 + 3 \zeta_{5}^3 +2 \zeta_{5}^4 
\end{eqnarray*}
and hence  $a_1 = 0$, $a_2 = -6$, $a_3 = 3$, $a_4 = 2$.
\vskip2mm \noindent
Substituting the values of $a_i$'s, we get
\begin{align*}\nonumber
	G = D^t &=\left(
	\begin{array}{cccc}
		9 & 1 & 0 & 7  \\
		2 & 3 & 55 & 0 
	\end{array}
	\right)
\end{align*}
a generating matrix of an MDS code of type $[4, 2, 3]$.
\begin{align*}\nonumber
	Y &=\left(
	\begin{array}{cc}
		9  & 	1  \\
		2	&  3
	\end{array}
	\right)
\end{align*}
Then  determinant  of $Y$ is $\triangle= 25 \neq 0$ and 
\begin{align*}\nonumber
	Y^{-1} &=\displaystyle\frac{1}{25}\left(
	\begin{array}{cc}
		3  &  60  \\
		59	& 9
	\end{array}
	\right) = \left(
	\begin{array}{cc}
		5  & 	39  \\
		17	& 15
	\end{array}\right)
\end{align*} 
We get a generator matrix in the standard form:
\begin{align*}\nonumber
	G'=Y^{-1}G &=\left(\begin{array}{cc} 
		I_2 & \left(
		\begin{array}{cc}
			5  & 	39  \\
			17	& 15
		\end{array}
		\right) \cdot\left(\begin{array}{cc}
			0 & 7 \\
			55 &  0 
		\end{array} \right)
	\end{array}	\right)=\left(\begin{array}{cc} 
	I_2 & \left(
	\begin{array}{cc}
		10  & 	35   \\
		32	& 58
	\end{array}
	\right)
	\end{array}\right)
\end{align*}
The parity check matrix is 
\begin{align*}\nonumber
	H &= \left(\begin{array}{cc} 
		\displaystyle\frac{1}{25} \cdot\left(
		\begin{array}{cc}
			0 & 55 \\
			7 &  0 
		\end{array}
		\right) \cdot \left(
		\begin{array}{cc}
			58  & 	2  \\
			1	& 52
		\end{array}
		\right) & I_2
	\end{array} \right)= \left(\begin{array}{cc} 
 	\left(\begin{array}{cc}
 		51 & 29 \\
 		26 &  3 
 	\end{array}
 	\right)  & I_2
 \end{array} \right)
\end{align*} 
Syndrome of $v=[a, b, c, d]$ is 
\begin{align*}\nonumber
	vH^t &=(a, b, c, d)\cdot \left(
	\begin{array}{cc}
		51 & 26 \\
		29 &  3 \\
		1 &0\\
		0& 1
	\end{array}
	\right)
\end{align*}
Here $A_1=51$, $A_2=26$, $A_3=29$, $A_4=3$.
Let $w=(11,4,55,7)$	be a codeword. We illustrate decoding using above MDS-code.
\begin{center}
	\begin{tabular}{c|c|c|c|c}
		Received word $v$ & syndrome of $v$ & $e_0$ & error vector $e$ & codeword $w=v-e$\\
		\hline
		$(9,4,55,7)$ & $(20,9)=59(51, 26)$& $59$ & $(59,0,0,0)$ & $(11,4,55,7)$\\
		$(11,17,55,7)$ & $(11,39)=13(29, 3)$& $13$ & $(0,13,0,0)$  & $(11,4,55,7)$\\
		$(11,4,19,7)$ & $(25,0)=25(1, 0)$& $25$ & $(0,0,25,0)$ & $(11,4,55,7)$\\
		$(11,4,55,18)$ & $(0,11)=11(0, 1)$& $11$ & $(0,0,0,11)$  & $(11,4,55,7)$
	\end{tabular}	
\end{center}

\end{example}
\vskip2mm
\noindent \textbf{Future Scope:} We have thus verified the conjecture in \S 3 for orders $3, 5$ and thereby obtained Gauss-Dickson codes. It is expected that these results can be carried forward for higher values of $l$, 
however the calculations become laborious even using a software. Such results are expected for Jacobi sums $J(1, n) $ 
whenever the Jacobi sums has distinct conjugates. Also Jacobi codes of composite order can be tried. Vikas Jadhav and Katre have observed
that for $p=79$ and $l =13$, we do not get MDS codes for certain generators of $\mathbb F_p^*$, however for other
$p \equiv 1 \Mod{13}$ we get MDS codes. Thus there is a possibility of exceptional primes for the conjecture in $\S 3$.

\vskip2mm \noindent
{\bf Acknowledgement}: The first author acknowledges support from Lokmanya Tilak Chair of  S. P. Pune University.  
The second  author thanks UGC-Teacher fellowship received through Nowrosjee Wadia College, Pune to work at BP, Pune, during which a part of this work was done. 
 
\end{document}